# Dimensionality reduction of networked systems with separable coupling-dynamics: theory and applications

**Authors:** Chengyi Tu[1,2*], Ying Fan[3], Tianyu Shi[2*]

[1]Department of Environmental Science, Policy, and Management, University of California, Berkeley; Berkeley, 94720, USA.

[2]School of Economics and Management, Zhejiang Sci-Tech University; Hangzhou, 310018, China.

[3]College of Geography and Environment, Shandong Normal University; Jinan, 250358, China.

*Corresponding author. Email: chengyitu1986@gmail.com and shitianyu@zstu.edu.cn

# Abstract

Complex dynamical systems are prevalent in various domains, but their analysis and prediction are hindered by their high dimensionality and nonlinearity. Dimensionality reduction techniques can simplify the system dynamics by reducing the number of variables, but most existing methods do not account for networked systems with separable coupling-dynamics, where the interaction between nodes can be decomposed into a function of the node state and a function of the neighbor state. Here, we present a novel dimensionality reduction framework that can effectively capture the global dynamics of these networks by projecting them onto a low-dimensional system. We derive the reduced system's equation and stability conditions, and propose an error metric to quantify the reduction accuracy. We demonstrate our framework on two examples of networked systems with separable coupling-dynamics: a modified susceptible-infected-susceptible model with direct infection and a modified Michaelis-Menten model with activation and inhibition. We conduct numerical experiments on synthetic and empirical networks to validate and evaluate our framework, and find a good agreement between the original and reduced systems. We also investigate the effects of different network structures and parameters on the system dynamics and the reduction error. Our framework offers a general and powerful tool for studying complex dynamical networks with separable coupling-dynamics.

**Keywords**: dimensionality reduction; complex dynamical networks; separable coupling-dynamics

# Introduction

Complex dynamical systems are ubiquitous in nature and society, permeating a wide array of fields, from biological networks to social systems [1-6]. These systems are integral to our understanding of the world around us, and their behavior and evolution present a fundamental challenge in the realms of science and engineering. However, the analysis and prediction of these systems are significantly complicated by their high dimensionality and nonlinearity [7, 8]. The sheer number of variables involved in these systems, coupled with the intricate ways in which they interact, can make it incredibly difficult to discern patterns or predict future behavior. This complexity is further compounded by the fact that these systems often exhibit nonlinearity, meaning that changes in input do not result in proportional changes in output. This lack of proportionality can lead to unexpected and chaotic behavior, making prediction and control even more challenging. Given these challenges, dimensionality reduction techniques have emerged as essential tools in the study of complex dynamical systems [9, 10]. These techniques aim to simplify the system dynamics by reducing the number of variables under consideration, thereby making the system more manageable and its behavior more understandable. By distilling the system down to its most critical components, researchers can gain insights into the underlying mechanisms that drive the system's behavior.

The dimensionality reduction of complex systems, despite significant efforts, continues to pose challenges, particularly in the context of complex networks [11, 12]. Various strategies, including linear projection, manifold learning, network embedding, and critical slowing down, have been explored to address these challenges. Gao et al. pioneered the application of dimensionality reduction to networked systems by proposing an analytical tool that leverages mean-field approaches to elucidate the effective parameters of high-dimensional networked systems [13]. However, Tu et al. identified a new condition that significantly restricts the applicability of Gao's framework [14]. Kundu et al. examined the precision of Gao's framework in the context of networks without degree correlation [15]. Subsequently, Tu et al. expanded Gao's framework to investigate the collapse or operation of networked systems, taking into account the complete heterogeneity in node-specific self- and coupling-dynamics, discrete-time dynamical systems and stochastic dynamical systems [16-18]. They also applied their framework to reduce a phenomenological plant-pollinator model with empirical networks to effective dynamics to reveal that the systems may undergo critical transitions to a state of low species abundance, which is undesirable for biodiversity and ecosystem functioning. To cater to a general network structure, Laurence et al. proposed a polynomial approximation framework for simplifying complex networks [19, 20]. Masuda et al. put forth a theory suggesting that a non-leading eigenvector of the adjacency matrix often yields better accuracy when applied to Laurence's framework [21]. Wu et al. proposed an enhanced method for dimension reduction, which leverages the concept of information entropy to forecast the resilience of networks [22]. Ghosh et al. developed a one-dimensional reduced model to analyze susceptible-infected-susceptible dynamics, incorporating higher-order interactions into their model [23].

However, most of these methods do not account for networked systems where the coupling-dynamics between nodes are separable, meaning that they can be decomposed into a function of the node state and a function of the neighbor state [24-26]. This characteristic allows for a more granular understanding of the system, as it enables the isolation and examination of individual node behaviors and their subsequent impact on their neighbors [27]. Despite the potential insights that could be gleaned from analyzing systems with separable coupling-dynamics, current methodologies largely overlook this aspect. This gap in the literature presents an opportunity for the development of novel analytical tools and frameworks that can harness the intricacies of separable coupling-dynamics to provide a more comprehensive understanding of networked systems. The exploration of such methodologies could have far-reaching implications across various fields, from understanding the spread of information or disease in social networks, to optimizing the performance of interconnected systems in engineering and technology. As such, there is a pressing need for further research in this area to unlock the full potential of network analysis in systems with separable coupling-dynamics.

Here, we present a novel dimensionality reduction framework that can effectively capture the global dynamics of these networks by projecting them onto a low-dimensional system. We derive the reduced system's equation and stability conditions, and propose an error metric to quantify the reduction accuracy. We demonstrate our framework on two examples of networked systems with separable coupling-dynamics: a modified susceptible-infected-susceptible model with direct infection and a modified Michaelis-Menten model with activation and inhibition. We conduct numerical experiments on synthetic and empirical networks to validate and evaluate our framework, and find a good agreement between the original and reduced systems. We also investigate the effects of different network structures and parameters on the system dynamics and the reduction error. Our framework offers a general and powerful tool for studying complex dynamical networks with separable coupling-dynamics.

# Methods

## Dimensionality reduction in separably networked dynamic systems

We study a class of networked systems with $N$ nodes whose states $\mathbf{x} = (x_1, \ldots, x_N)^T$ evolve according to the equation

$$\frac{dx_i}{dt} = F_i(x_i) + \sum_{j}^{N} A_{ij} G_{ij}(x_i, x_j) \quad (1)$$

where $F_i(x_i)$ self-dynamics function and $G_{ij}(x_i, x_j)$ is the coupling function between node $i$ and its neighbors $j$, as determined by the adjacency matrix $\mathbf{A} \in R^{N \times N}$, which encodes the network structure with $A_{ij}$ strength of the connection from node $j$ to node $i$. In the context of this study, we postulate that the coupling function, denoted as $G_{ij}(x_i, x_j)$, exhibits separability. This implies that it can be expressed as the product of two distinct functions, namely $P(x_i)$ and $Q(x_j)$, each solely dependent on the state of nodes $i$ and $j$ respectively. Consequently, the coupling function can be represented as $G_{ij}(x_i, x_j) = P_i(x_i) * Q_j(x_j)$. This property simplifies the analysis of the network dynamics and enables us to explore the effects of different coupling schemes [24, 25, 27].

We introduce an operator $\mathcal{L}(\mathbf{x}) = \frac{1}{N}\sum_{j=1}^{N} s_j^{out} x_j / \frac{1}{N}\sum_{j=1}^{N} s_j^{out} = \frac{\langle \mathbf{s}^{out} \cdot \mathbf{x} \rangle}{\langle \mathbf{s}^{out} \rangle}$, where $\mathbf{s}^{out} = (s_1^{out}, \ldots, s_N^{out})$ is the vector of the out-degree of the network $\mathbf{A}$ [13]. The operator $\mathcal{L}$ is compatible with both linear time-invariance (LTI) functions and the Hadamard product, i.e., $\mathcal{L}(a\mathbf{x} + b\mathbf{y}) = a\mathcal{L}(\mathbf{x}) + b\mathcal{L}(\mathbf{y})$ and $\mathcal{L}(\mathbf{x} \circ \mathbf{y}) \approx \mathcal{L}(\mathbf{x})\mathcal{L}(\mathbf{y})$ [16]. If the degree correlation of network $\mathbf{A}$ is negligible (the neighborhood of node $i$ mirrors the neighborhoods of all other nodes), then the expression

$$\sum_{j}^{N} A_{ij} G_{ij}(x_i, x_j) \approx s_i^{in} * \mathcal{L}(P_i(x_i) * \mathbf{Q}(\mathbf{x})) \approx s_i^{in} * P_i(x_i) * \mathcal{L}(\mathbf{Q}(\mathbf{x}))$$ where $\mathbf{Q}(\mathbf{x}) = (Q_1(x_1) \ldots Q_N(x_N))$.

The simplification of the term $\mathcal{L}(\mathbf{Q}(\mathbf{x}))$ can be achieved by expressing each $Q_j(x_j)$ as a linear combination of $m$ subfunctions, i.e., $Q_j(x_j) = b_{j,1}q_1(x_i) + b_{j,2}q_2(x_i) + \cdots + b_{j,m}q_m(x_i)$. By applying the operator $\mathcal{L}$, which is linear time-invariance, we obtain $\mathcal{L}(\mathbf{Q}(\mathbf{x})) = \mathcal{L}\begin{pmatrix} Q_1(x_1) \\ \vdots \\ Q_N(x_N) \end{pmatrix} = \mathcal{L}\begin{pmatrix} b_{1,1}q_1(x_1) \\ \vdots \\ b_{N,1}q_1(x_N) \end{pmatrix} + \ldots + \mathcal{L}\begin{pmatrix} b_{1,m}q_m(x_1) \\ \vdots \\ b_{N,m}q_m(x_N) \end{pmatrix}$. Then, using the operator $\mathcal{L}$ which exhibits the Hadamard product property, we can derive an approximation of

$$\mathcal{L}(\mathbf{Q}(\mathbf{x})) = \mathcal{L}(B^1 \circ q_1(\mathbf{x})) + \cdots + \mathcal{L}(B^m \circ q_m(\mathbf{x})) \approx \mathcal{L}(B^1)\mathcal{L}(q_1(\mathbf{x})) + \cdots + \mathcal{L}(B^m)\mathcal{L}(q_m(\mathbf{x}))$$ where $B^k = (b_{1,k}, \ldots, b_{N,k})^T$ is the $k$-th column of matrix $\mathbf{B}$. Given that the subfunctions are identical across all nodes, an approximation of can be further derived

$$\mathcal{L}(\mathbf{Q}(\mathbf{x})) \approx \mathcal{L}(B^1)q_1(\mathcal{L}(\mathbf{x})) + \ldots + \mathcal{L}(B^m)q_m(\mathcal{L}(\mathbf{x})) \approx \mathcal{L}(B^1)q_1(x_{eff}) + \ldots + \mathcal{L}(B^m)q_m(x_{eff}) = \sum_{k=1}^{m} B_{eff}^k q_k(x_{eff})$$ . In

instances where the function $Q_i(x_i)$ may not be a linear combination of subfunctions or may vary across nodes, Chebyshev polynomials can be employed to approximate $Q_i(x_i)$ as $\sum_{k=1}^{m} b_{i,k} x^{(k-1)}$, thereby minimizing the approximation error [28, 29]. The operator $\mathcal{L}$ can then be applied to obtain $\mathcal{L}(\mathbf{Q}(\mathbf{x})) = \sum_{k=1}^{m} B_{eff}^k x_{eff}^{(k-1)}$, which we designate as $Q_{eff}$. In a similar vein, the aforementioned process can be employed to approximate $\mathcal{L}(\mathbf{P}(\mathbf{x}))$ and $\mathcal{L}(\mathbf{F}(\mathbf{x}))$, which are denoted as $P_{eff}$ and $F_{eff}$, respectively.

We can rewrite Eq. (1) in the following form $\frac{dx_i}{dt} \approx F_i(x_i) + s_i^{in} * P_i(x_i) * Q_{eff}$, with its vector notation being $\frac{d\mathbf{x}}{dt} = \mathbf{F}(\mathbf{x}) + P_{eff} * \mathbf{s}^{in} \circ \mathbf{Q}(\mathbf{x})$. By applying the operator $\mathcal{L}$ to both sides of the vector equation, we obtain the effective equation of Eq. (1):

$$O = \frac{dx_{eff}}{dt} \approx F_{eff} + A_{eff} P_{eff} Q_{eff} \quad (2)$$

## Stability of solution of effective equation

The stability of the solutions derived from the effective equation $x_{eff}^*$ can be associated with a region defined by the following set of equations:

$$\begin{cases} O(x_{eff}^*) = 0 \\ Re\left[\frac{\partial O}{\partial x_{eff}}\bigg|_{x_{eff}=x_{eff}^*}\right] < 0 \end{cases} \quad (3)$$

By applying the analytical method to the low-dimensional effective system given by Eq. (3), we can derive the resilience function that characterizes the possible states of the system in terms of the effective parameters. This function enables us to examine the stability of the system and the occurrence of critical transitions, which are abrupt changes in the system's state caused by small perturbations. To test the validity of this method, we compare the results obtained from the low-dimensional system with those from the high-dimensional system described by Eq. (1), which captures the full complexity of the system dynamics.

## Error definition

We evaluate the performance of our framework by computing an error metric that measures the projection distance between the numerical and theoretical solutions. The numerical solution $x_{eff}^N$ is obtained by the original high-dimensional system, Eq. (1), which gives the state vector $\mathbf{x}$ of the system and definition $x_{eff} = \frac{\langle \mathbf{s}^{out} . \mathbf{x} \rangle}{\langle \mathbf{s}^{out} \rangle}$. The theoretical solution $x_{eff}^A$ is derived from the low-dimensional effective equation, Eq. (3). The projection distance from point obtained by numerical solution to the surface obtained by theoretical solution can be defined as an error

$$err = \left| \frac{x_{eff}^N - x_{eff}^A}{x_{eff}^N} \right| \quad (4)$$

The framework relies on the assumption that the numerical simulations can approximate the theoretical predictions with a negligible error. This assumption is valid when the error is sufficiently small, such that the simulated point is close to the true point in the parameter space. However, when the error is large, the framework may fail to converge to the right solution or produce inaccurate estimates.

# Results

## Modified susceptible-infected-susceptible model with random networks

The susceptible-infected-susceptible (SIS) model is a widely used framework for studying the dynamics of infectious diseases that do not confer immunity upon recovery, such as the common cold or sexually transmitted infections [30, 31]. The model assumes that individuals in a population can be classified into two states: susceptible (S) or infected (I). Susceptible individuals can become infected when they have contact with infected individuals, and infected individuals can recover and return to the susceptible state. The SIS model is an endemic model, meaning that the disease persists in the population indefinitely. However, the SIS model does not account for the possibility that some individuals may become infected by external sources, such as environmental factors or imported cases [32, 33]. To capture this effect, we propose a modified SIS model that incorporates a term for direct infection, which represents the probability that a susceptible individual becomes infected without contact with an infected individual. The modified SIS model is described by the following equation:

$$\frac{dx_i}{dt} = -e_i x_i + \sum_{j}^{N} A_{ij}(1-x_i)x_j^Z \quad (5)$$

where $N$ is network size or individual number, $0 \leq x_i \leq 1$ denotes the probability that node $i$ is infected at a given time, $e_i$ is the rate at which node $i$ recovers from the infection and becomes susceptible again, matrix $\mathbf{A}$ captures the infection rates between nodes where $A_{ij}$ is the rate at which node $i$ becomes infected by node $j$ if they are connected, and $Z \sim B(p)$ is a Bernoulli distribution gives value $z=1$ with probability $p$ and $z=0$ with probability $1-p$ where $p$ is the probability of direct infection, which means that a susceptible individual become infected by contacting with an infected individual. The first term on the right-hand side of Eq. (5) accounts for the process of recovery, and the second term accounts for the process of infection. The infection process depends on the parameter of Bernoulli distribution and the state of the neighboring nodes. If $z = 0$, then node $i$ is directly infected regardless of the state of its neighbors. If $z = 1$, then node $i$ is infected only if it has an infected neighbor $j$. Therefore, we define the parameter $p$ in the Bernoulli distribution as the direct infection rate, which measures the likelihood that a susceptible individual becomes infected without contact with an infected individual. The effective equation of the modified SIS model is obtained by taking the expectation of Eq. (5) over the direct infection rate $p$, which yields:

$$\frac{dx_{eff}}{dt} = -e_{eff} x_{eff} + A_{eff}(1 - x_{eff}) x_{eff}^{p} \quad (6)$$

We perform numerical simulations to test the validity of the effective equation derived from the modified SIS model with non-direct infection. We consider a network of $N$ nodes, where each node had a random recovery rate $e_i$ drawn from a uniform distribution between 0 and $2\mu_e$, with $\mu_e$ being the mean recovery rate. We set $\mu_e = 100$ and $\sigma_e = |\mu_e / 3|$ as the default values. We also vary the probability of direct infection $p$ from 0.25 to 0.75, and the network size $S$ from 100 to 300. We construct some Erdős–Rényi (ER) networks with a connectivity of $C = 0.5$ as the underlying network structure. We initialize the state of each node $x_i$ with a random value between 0 and 0.1 for a low initial population, or between 0.9 and 1 for a large initial population. We compare the numerical solutions of the modified SIS model (Eq. (5)) and the analytical solution of effective equation (Eq. (6)) for different values of $p$ and $S$. The results are shown in Fig. 1.

The numerical simulations confirm that the effective equation accurately capture the dynamics of the modified SIS model with non-direct infection for the given parameter configurations, network structures, and initial conditions. The effective equation predicted that the system had a unique stable equilibrium, which was consistent with the numerical solutions of the modified SIS model (see Fig. 1 a-c). The error between the numerical solutions and the effective equation is small, especially for small values of $p$ and large values of $S$ (see Fig. 1 d-f). The error decreases as $p$ decreases or $S$ increases, indicating that the effective equation was more accurate for higher probabilities of direct infection or larger network sizes. The error approaches zero as $S$ tends to infinity, suggesting that the effective equation was asymptotically exact for infinite networks.

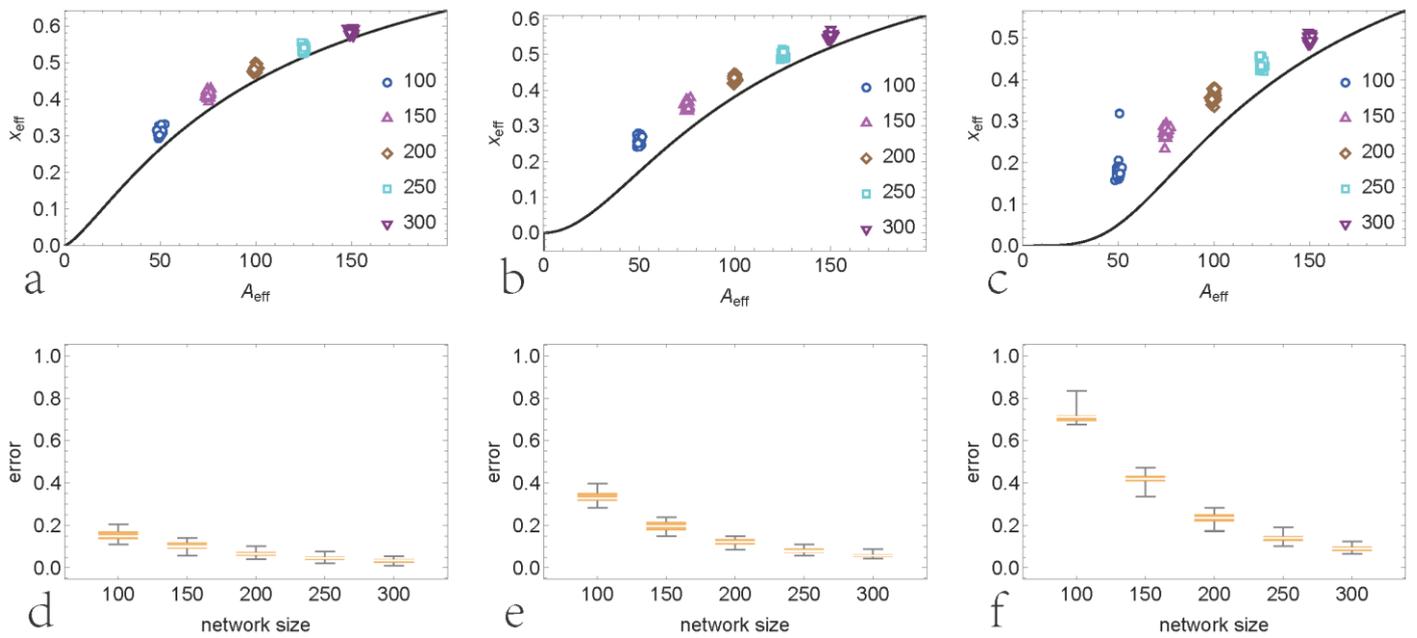

**Fig. 1 | Comparison of numerical and theoretical results of modified SIS model with non-direct infection on ER networks with varying network sizes and probabilities of direct infection.** (a-c) The probability of direct infection is 0.25, 0.5 and 0.75, respectively. The solid black line represents the analytical prediction, and the colored markers represent the numerical simulations for each network size. Each case is carried out 50 simulations. (d-f) The box-and-whisker chart of error corresponding to a-c.

We extend our analysis to investigate the influence of different social network structures on the dynamics of the modified SIS model with non-direct infection. We assume that the nodes are connected by an Erdős–Rényi (ER) network, where each pair of nodes has an equal probability of being connected, as a baseline scenario. We explore two alternative scenarios, where the nodes are connected by a Barabási–Albert (BA) network or a small-world (SW) network, which reflect some of the properties of real-world social networks. The BA network has a scale-free degree distribution, meaning that some nodes have a very high number of connections, while most nodes have a very low number of connections. This structure arises from a preferential attachment mechanism, where new nodes tend to connect to existing nodes with higher degrees. The SW network has a high clustering coefficient and a low average path length, meaning that the nodes are locally connected and globally reachable. This structure mimics the features of many natural and social systems, where local interactions and global connectivity coexist.

We perform numerical simulations to test the validity of the effective equation derived from the modified SIS model with non-direct infection for different network structures. The parameter configurations, network structures, and initial conditions are the same as ER case. The results are shown in Fig. 2 and 3. The results confirm that the effective equation accurately capture the dynamics of the modified SIS model for the given parameter configurations, network structures, and initial conditions. We find consistent results between the numerical simulation and the theoretical prediction, so our main conclusions are robust and do not depend significantly on the type of social network connecting the nodes.

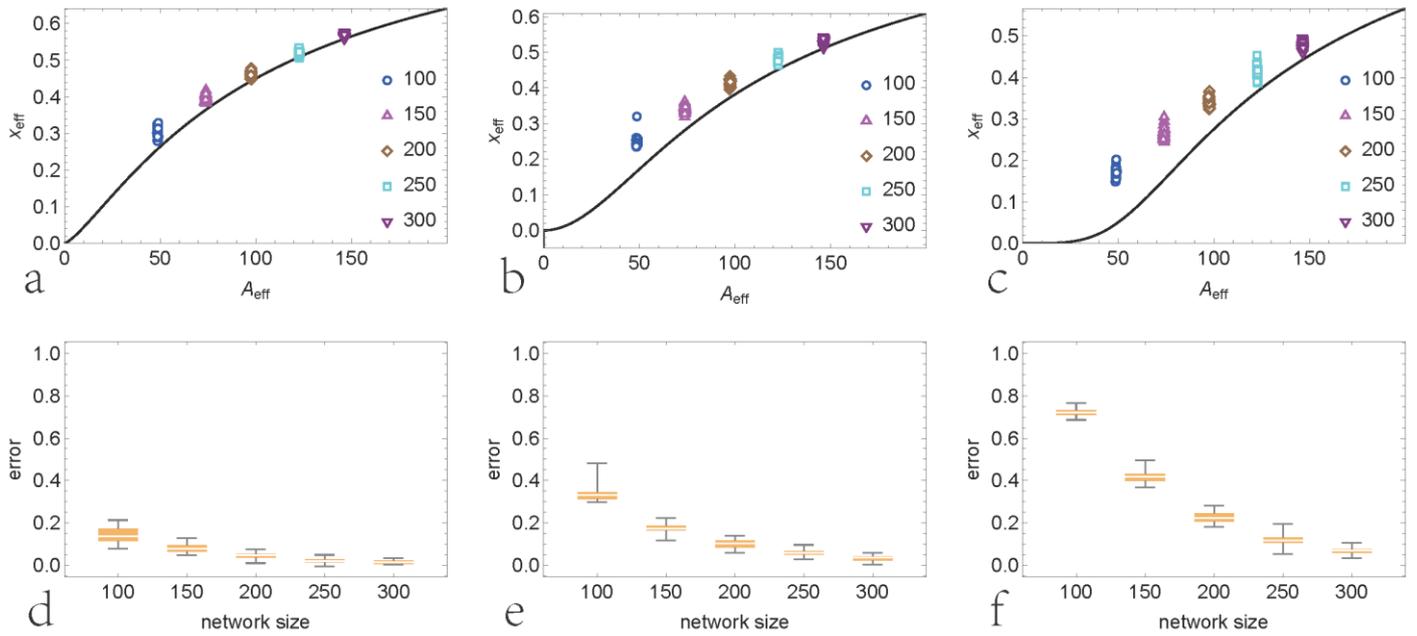

**Fig. 2 | Comparison of numerical and theoretical results of modified SIS model with non-direct infection on BA networks with varying network sizes and probabilities of direct infection.** (a-c) The probability of direct infection is 0.25, 0.5 and 0.75, respectively. The solid black line represents the analytical prediction, and the colored markers represent the numerical simulations for each network size. Each case is carried out 50 simulations. (d-f) The box-and-whisker chart of error corresponding to a-c.

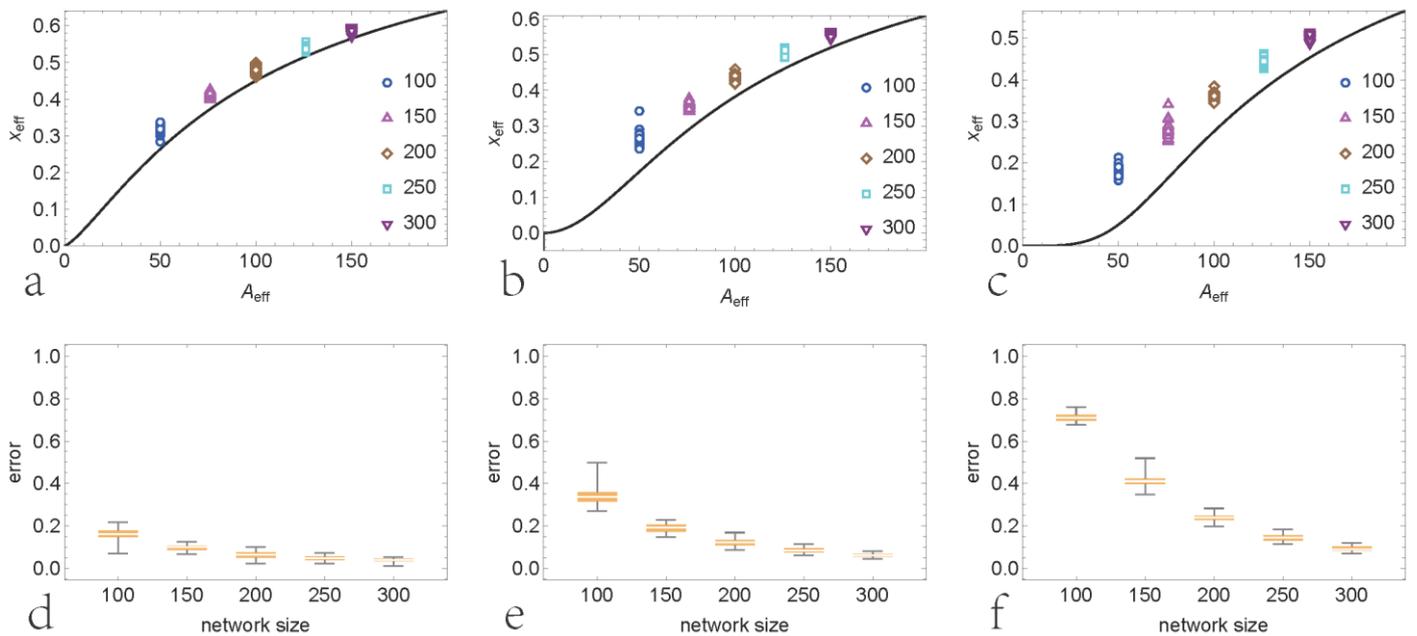

**Fig. 3 | Comparison of numerical and theoretical results of modified SIS model with non-direct infection on SW networks with varying network sizes and probabilities of direct infection.** (a-c) The probability of direct infection is 0.25, 0.5 and 0.75, respectively. The solid black line represents the analytical prediction, and the colored markers represent the numerical simulations for each network size. Each case is carried out 50 simulations. (d-f) The box-and-whisker chart of error

corresponding to a-c.

## Modified susceptible-infected-susceptible model with empirical networks

We also extend our analysis to examine the dynamics of the modified SIS model with non-direct infection on a real-world contact network. We use the infect-dublin dataset, which is a human contact network where nodes represent humans and edges represent proximity in the physical world [34, 35]. The network contains 410 nodes and 2765 edges. We use the same parameter configurations and initial conditions as in the previous scenarios. We compare the numerical solutions of the modified SIS model and the effective equation for the real-world network. The results are shown in Fig. 4.

Our results demonstrate that the effective equation provides a reliable approximation of the modified SIS model with non-direct infection for the real-world network. The numerical solutions and the effective equation agree with each other, and they show similar trends of infection prevalence and error. The error diminishes as the probability of direct infection decreases, consistent with the previous scenarios. Therefore, we infer that our main results are robust and do not depend significantly on the type of network connecting the individuals.

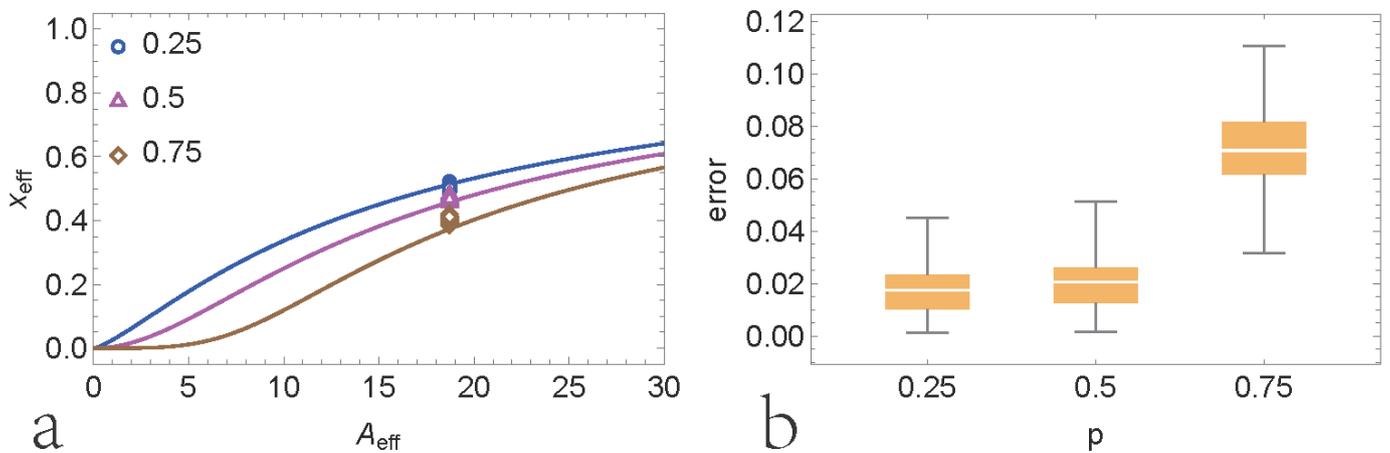

**Fig. 4 | Comparison of numerical and theoretical results of modified SIS model with non-direct infection on an empirical network with vary probabilities of direct infection.** (a-c) The probability of direct infection is 0.25, 0.5 and 0.75, respectively. The colored line represents the analytical prediction, and the colored markers represent the numerical simulations for each network size. Each case is carried out 50 simulations. (d-f) The box-and-whisker chart of error corresponding to a-c.

## Modified Michaelis-Menten model with empirical networks

To illustrate the applicability of our method to biological systems, we consider a gene regulatory network modelled by Michaelis-Menten (MM) kinetics [36-38]. The MM kinetics are widely used to describe the rate of enzymatic reactions as a function of the substrate concentration, and can capture the nonlinear effects of gene regulation, such as activation and inhibition. However, most existing works assume homogeneous coupling-dynamics between nodes, meaning that the interaction type is the same for all connected gene $j$ for each gene $i$. This assumption neglects the heterogeneity of gene interactions, which can be either activating or inhibiting. To account for this heterogeneity, we propose a modified MM model that incorporates a term for activation and inhibition, which quantifies the probability that two genes are activated or inhibited by each other. The modified MM model is given by the following equation:

$$\frac{dx_i}{dt} = -e_i x_i + \sum_j^N A_{ij} \frac{x_j^Z}{x_j + 1} \quad (7)$$

where $N$ is network size or gene number, $x_i$ is the expression level of gene $i$, $e_i$ is the self-dynamic rate, $\mathbf{A}$ is the adjacency matrix that represents the network topology, and $Z \sim B(p)$ is a Bernoulli distribution gives value $z = 1$ with probability $p$ and $z = 0$ with probability $1-p$ where $p$ is the probability of activation, which means that the gene $i$ is activated by gene $j$. The first term is self-regulation function of gene $i$, the second term is the interaction between gene $j$ and gene $i$ which is modeled by Hill functions [39, 40]. The Hill functions are commonly used to describe the sigmoidal response of biochemical systems: If $z = 0$, then the coupling dynamics become $\frac{1}{x_j + 1}$, where means that gene $i$ is inhibited by gene $j$. If $z = 1$, then the coupling dynamics become $\frac{x_j}{x_j + 1}$, where means that gene $i$ is activated by gene $j$. Therefore, our model can capture the heterogeneity of gene interactions, which can be either activating or inhibiting. The effective equation of the modified MM model is obtained by taking the expectation of Eq. (7) over the activation rate $p$, which yields:

$$\frac{dx_{eff}}{dt} = -e_{eff} x_{eff} + A_{eff} \frac{x_{eff}^p}{1 + x_{eff}} \quad (8)$$

We perform numerical simulations to test the validity of the effective equation derived from the modified MM model with inhibition and activation. We consider a network of $N$ nodes, where each node had a random self-dynamic rate $e_i$ drawn from a uniform distribution between 0 and $2\mu_e$, with $\mu_e$ being the mean recovery rate. We set $\mu_e = 8$ and $\sigma_e = |\mu_e / 3|$ as the default values. We also vary the probability of activation rate $p$ from 0.25 to 0.75. We use bio-diseasome, a network of human diseases and their associated genes, based on the concept of the human diseasome and the human disease network [41, 42]. The bio-diseasome network contains 516 nodes representing diseases and 1188 edges representing gene-disease associations. The network can be used to explore the relationships between diseases and genes, as well as the properties and statistics of the network structure. We initialize the state of each node $x_i$ with a random value between 0 and 0.1 for a low initial population, or between 0.9 and 1 for a large initial population. We compare the numerical solutions of the modified MM model (Eq. (7)) and the analytical solution of effective equation (Eq. (8)) for different values of $p$. The results are shown in Fig. 5.

The results show that the effective equation accurately captures the dynamics of the modified MM model for the real-world network. The numerical solutions and the effective equation are consistent with each other, and they exhibit similar patterns and error (see Fig. 5 a-c). The errors are nearly identical for different $p$ values, indicating that the effective equation is robust to the changes in the activation rate (see Fig. 5 d-f). These results demonstrate that the effective equation is a valid approximation of the modified MM model for the real-world network, and it can be used to analyze the network dynamics more efficiently.

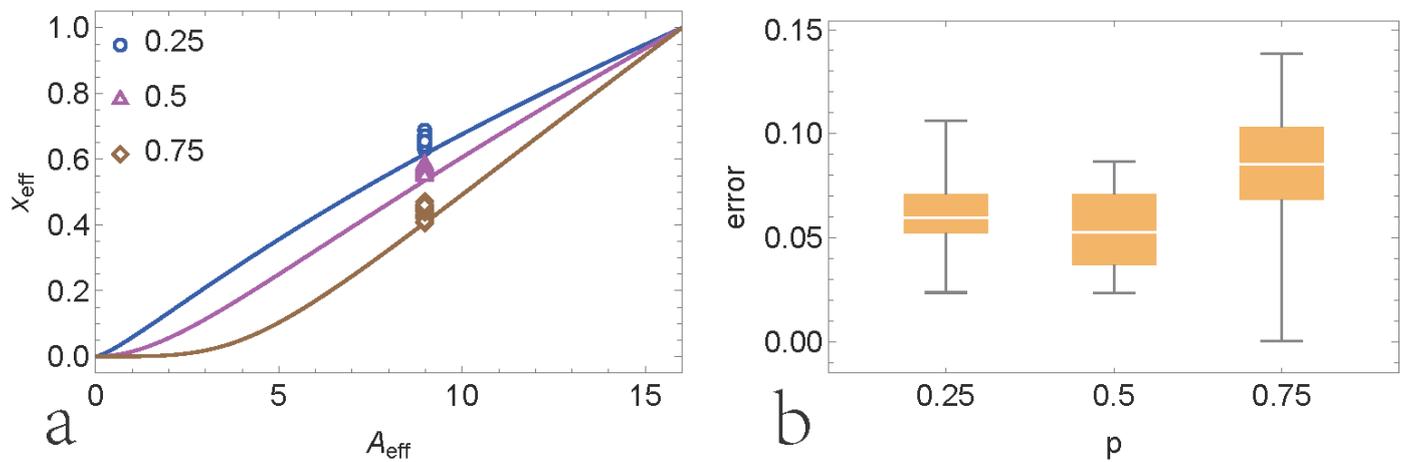

**Fig. 5 | Comparison of numerical and theoretical results of modified MM dynamics on an empirical network with varying probabilities of activation.** (a-c) The probability of activation is 0.25, 0.5 and 0.75, respectively. The colored line represents the analytical prediction, and the colored markers represent the numerical simulations for each network size. Each case is carried out 50 simulations. (d-f) The box-and-whisker chart of error corresponding to a-c.

# Discussions and conclusions

We present an analytical framework for dimensionality reduction of complex dynamical systems that are networked and have separable coupling-dynamics. The framework can effectively capture the global dynamics of these systems by projecting them onto a low-dimensional system, which can be derived from the original system using a set of operators and approximations. The framework can handle various types of systems, such as biological, ecological, and social systems, and can account for different network structures and parameters. We claim that the effective equation is a valid approximation of the original system, and that it can be used to analyze the system dynamics more efficiently and accurately.

We demonstrate the applicability and validity of the method by applying it to two examples: a modified susceptible-infected-susceptible (SIS) model with non-direct infection, and a modified Michaelis-Menten (MM) model with both activation and inhibition. The modified SIS model incorporates a term for non-direct infection, which represents the probability that a susceptible individual becomes infected without contact with an infected individual. The modified MM model incorporates a term for activation and inhibition, which represents the possibility that two genes are activated or inhibited by each other. We compare the numerical solutions of the original high-dimensional systems with the analytical solutions of the effective equations for different parameter configurations, network structures, and initial conditions. The results show that the effective equations accurately capture the dynamics of the original systems, and that the error between the numerical and analytical solutions is small and decreases as the network size increases.

The paper makes a significant contribution to the field of complex systems, as it provides a general and powerful framework for dimensionality reduction of networked dynamical systems with separable coupling-dynamics. The method can be applied to a wide range of systems, such as epidemic, ecological, social, and biological systems, and can help to reveal the universal patterns and mechanisms underlying the system dynamics. The method can also facilitate the analysis of the system resilience and the detection of the critical transitions, which are important for understanding and managing the system behavior.

However, the paper also has some limitations and open questions that need to be addressed in future work. First, the paper assumes that the coupling-dynamics are random and independent of the local dynamics, which may not be realistic for some systems. For example, in the modified SIS model, the non-direct infection rate may depend on the environmental factors or the imported cases, which may vary over time and space. In the modified MM model, the activation and inhibition rates may depend on the gene expression levels or the regulatory feedbacks, which may be nonlinear and complex. Therefore, it would be interesting to explore how the method can be extended to account for the dependence and correlation of the coupling-dynamics. Second, the paper assumes that the coupling-dynamics are separable, which means that the interaction between two nodes only depends on their own states, and not on the states of other nodes. This assumption may not hold for some systems, such as neural networks, where the interaction between two neurons may depend on the activity of other neurons in the network. Therefore, it would be useful to investigate how the method can be generalized to handle non-separable coupling-dynamics. Third, the paper only considers two examples of networked dynamical systems, which may not be representative of the diversity and complexity of the real-world systems. Therefore, it would be desirable to apply the method to more examples of networked dynamical systems, such as climate, economic, and technological systems, and to compare the results with other existing methods for dimensionality reduction, such as principal component analysis, singular value decomposition, or autoencoders.

# Data and code availability

The ready-to-use notebook codes to reproduce the results presented in the current study are available in OSF with the access code q2n74 (https://osf.io/q2n74/).

# Acknowledgements

This work was supported by Zhejiang Provincial Natural Science Foundation of China (Grant No. LZ22G010001), Natural Science Fund of Zhejiang Province, China (LQ19G010004), and Science Foundation of Zhejiang Sci-Tech University (ZSTU) under Grant No. 18092125-Y and 22092034-Y.